\documentclass[12pt]{article}
\usepackage[pctex32]{graphics}
\usepackage{amsthm,amsmath,amssymb,amscd,verbatim,epsfig}    
\usepackage{amssymb}
\usepackage{graphicx}
\newcommand{\beq}{\begin{equation}}
\newcommand{\eeq}{\end{equation}}

\date{}

\newcommand{\f}{\frac}
\newcommand{\ra}{\rightarrow}

\newcommand{\sub}{\subset}
\newcommand{\supp}{\supset}

\newcommand{\sq}{$\blacksquare$}

\begin{document}
\title{A\ collection\ of\ simple\ proofs\ of\ Sharkovsky's\ theorem}
\author{Bau-Sen Du \\ [.2cm]
Institute of Mathematics \\
Academia Sinica \\
Taipei 11529, Taiwan \\
dubs@math.sinica.edu.tw \\}
\maketitle
\begin{abstract}
Based on various strategies, we obtain several simple proofs of the celebrated Sharkovsky cycle coexistence theorem. 
\end{abstract}

\section{Introduction}
Throughout this note, $I$ is a compact interval, and $f : I \ra I$ is a continuous map.   For each integer $n \ge 1$, let $f^n$ be defined by: $f^1 = f$ and $f^n = f \circ f^{n-1}$ when $n \ge 2$.  For $x_0$ in $I$, we call $x_0$ a periodic point of $f$ with least period $m$ (or a period-$m$ point of $f$) if $f^m(x_0) = x_0$ and $f^i(x_0) \ne x_0$ when $0 < i < m$.  If $f(x_0) = x_0$, then we call $x_0$ a fixed point of $f$.  

For discrete dynamical systems defined by iterated interval maps, one of the most remarkable results is Sharkovsky's theorem 
which states as follows:

\noindent
{\bf Theorem (Sharkovsky{\bf{\cite{mi, sh}}}).}  {\it Let the Sharkovsky's ordering of the natural numbers be defined as follows: $$3 \prec 5 \prec 7 \prec  9 \prec \cdots \prec 2 \cdot 3 \prec 2 \cdot 5 \prec 2 \cdot 7 \prec  2 \cdot 9 \prec \cdots \prec 2^2 \cdot 3 \prec 2^2 \cdot 5 \prec 2^2 \cdot 7 \prec  2^2 \cdot 9 \prec \cdots$$ $$\prec \cdots \prec 2^3 \prec 2^2 \prec 2 \prec 1.$$  Then the following three statements hold:
\begin{itemize}
\item[\rm{(1)}] 
Assume that $f : I \to I$ is a continuous map.  If $f$ has a period-$m$ point and if $m \prec n$, then $f$ also has a period-$n$ point.

\item[\rm{(2)}]
For each positive integer $n$ there exists a continuous map $f : I \to I$ that has a period-$n$ point but has no period-$m$ point for any $m$ with $m \prec n$.

\item[\rm{(3)}]
There exists a continuous map $f : I \to I$ that has a period-$2^i$ point for $i = 0, 1, 2, \ldots$ but has no periodic point of any other period.
\end{itemize}}

There have been a number of different proofs {\bf{[1-5, 7-14, 16-31]}} of it in the past 30 years providing various viewpoints of this beautiful theorem.  Especially, some of these proofs are excellent for the discrete dynamical systems course.  So, why do we need another proof?  Sharkovsky's theorem is well-known for its simplicity in assumptions and yet abundance in conclusions.  Furthermore, what makes it more appealing is that its proof uses only the intermediate value theorem (in an indirect way through, say, directed graph arguments).  The idea was clear except that the details were a little messy {\bf{\cite{mi}}}.  Even the seemingly simplest case that if $f$ has a period-$m$ point with $m \ge 3$ then $f$ has a period-2 point cannot be explained in a few words {\bf{\cite{al, ba, c, str}}}.  Our aim is to find a proof of this theorem which uses the intermediate value theorem in a more intuitive and direct way (see section 3) so that this beautiful result can be someday introduced in the calculus curriculum.

It is well-known {\bf{\cite{du2, du3, str}}} that (1) is equivalent to the following three statements: 
\begin{itemize}
\item[(a)] if $f$ has a period-$m$ point with $m \ge 3$, then $f$ has a period-2 point; 

\item[(b)] if $f$ has a period-$m$ point with $m \ge 3$ and odd, then $f$ has a period-$(m+2)$ point; and 

\item[(c)] if $f$ has a period-$m$ point with $m \ge 3$ and odd, then $f$ has a period-$6$ point and a period-$(2m)$ point.  
\end{itemize}

\noindent
These three statements explain clearly why Sharkovsky's ordering is defined as it is.  Indeed, loosely speaking, by (b), we  have $3 \prec 5 \prec 7 \prec  9 \prec \cdots$, and then by (c), we obtain $3 \prec 5 \prec 7 \prec  9 \prec \cdots \prec 2 \cdot 3$.  Now with the help of Lemma 3 below, we obtain $3 \prec 5 \prec 7 \prec  9 \prec \cdots \prec 2 \cdot 3 \prec 2 \cdot 5 \prec 2 \cdot 7 \prec  2 \cdot 9 \prec \cdots \prec 2^2 \cdot 3$ and then the induction will prove the rest except the powers of 2 which follows from (a).  See section 13 for details.  Note that in (c) we don't require the existence of periodic points of all {\it even} periods.  Only the existence of period-6 and period-$(2m)$ points suffices.  This weaker requirement sometimes makes the arguments simpler, see sections 11 and 12.  

In this note, we present several strategies on how to prove (a), (b) and (c).  In section 2, we state some preliminary results.  In section 3, the proof is based on the existence of a point $v$ such that $f^2(v) < v < f(v)$ and $f(v)$ is in the given period-$m$ orbit.  In sections 4 and 5, we do some surgery to the map $f$ and then use Lemma 6 to prove (c).  In section 6, we examine how the iterates of a given period-$m$ point jump around a fixed point of $f$.  In section 7, we investigate how the points in a given period-$m$ orbit which lie on one side of a fixed point of $f$ correspond under the action of $f$.  In sections 8 and 9, we use the strategies in sections 6 and 7 to treat the special cases when $f$ has a periodic point of odd period $m \ge 3$ but no periodic points of odd periods $\ell$ with $\ell$ odd and $1 < \ell < m$.  In sections 10, 11, and 12, we concentrate our attention on the point $\min P$ and/or the point $\max P$ in any given period-$m$ orbit $P$.  The proof in each section is independent of the other.  Some proofs of (a), (b), or (c) in these sections can be combined to give various complete proofs of (a), (b) and (c).  For completeness, we also include here a slight modification of those proofs of Sharkovsky's theorem in {\bf{\cite{du2, du3}}}.

\section{Preliminary results}
To make this paper self-contained , we include the following results and some of their proofs.

\noindent
{\bf Lemma 1.}
If $f^m(x_0) = x_0$, then the least period of $x_0$ under $f$ divides $m$.

\noindent
{\bf Lemma 2.}
If $J$ is a closed subinterval of $I$ with $f(J) \supset J$, then $f$ has a fixed point in $J$.

\noindent
{\bf Lemma 3.}
{\it Let $k, m, n$, and $s$ be positive integers.  Then the following statements hold:
\begin{itemize}
\item[\rm{(1)}]
If $y$ is a periodic point of $f$ with least period $m$, then it is a periodic point of $f^n$ with least period $m/(m, n)$, where $(m, n)$ is the greatest common divisor of $m$ and $n$. 

\item[\rm{(2)}]
If $y$ is a periodic point of $f^n$ with least period $k$, then it is a periodic point of $f$ with least period $kn/s$, where $s$ divides $n$ and is relatively prime to $k$. 
\end{itemize}}

\noindent
{\it Proof.} (1) Let $t$ denote the least period of $x_0$ under $f^n$.  Then $m$ divides $nt$ since $x_0 = (f^n)^t(x_0) = f^{nt}(x_0)$.  Consequently, $\f m{(m, n)}$ divides $\f n{(m, n)} \cdot t$.  Since $\f m{(m, n)}$ and $\f n{(m, n)}$ are coprime, $\f m{(m, n)}$ divides $t$.  On the other hand, $(f^n)^{(m/(m, n))}(x_0) = (f^m)^{(n/(m, n))}(x_0) = x_0$.  So, $t$ divides $\f m{(m, n)}$.  This shows that $t = \f m{(m, n)}$.

(2) Since $x_0 = (f^n)^k(x_0) = f^{kn}(x_0)$, the least peirod of $x_0$ under $f$ is $\f {kn}s$ for some positive integer $s$.  By (1), $(\f {kn}s)/((\f {kn}s), n) = k$.  So, $\f ns = ((\f ns)k, n) = ((\f ns)k, (\f ns)s) = (\f ns)(k, s)$.  This shows that $s$ divides $n$ and $(s, k) = 1$.
\hfill\sq

\noindent
{\bf Lemma 4.}
Let $J$ and $L$ be closed subintervals of $I$ with $f(J) \supset L$.  Then there exists a closed subinterval $K$ of $J$ such that $f(K) = L$.

\noindent
{\it Proof.}
Let $L = [a, b]$.  Then since $\{ a, b \} \subset L \subset f(J)$, there are two points $p$ and $q$ in $J$ such that $f(p) = a$ and $f(q) = b$.  If $p < q$, let $c = \max \{ p \le x \le q : f(x) = a \}$ and let $d = \min \{ c \le x \le q : f(x) = b \}$.  If $p > q$, let $c = \max \{ q \le x \le p \mid f(x) = b \}$ and let $d = \min \{ c \le x \le p : f(x) = a \}$.  In any case, let $K = [c, d]$.  It is clear that $f(K) = L$.
\hfill\sq

If there are closed subintervals $J_0$, $J_1$, $\cdots$, $J_{n-1}, J_n$ of $I$ with $J_n = J_0$ such that $f(J_i) \supp J_{i+1}$ for $i = 0, 1, \cdots, n-1$, then we say that $J_0J_1 \cdots J_{n-1}J_0$ is a {\it{cycle of length}} $n$.  We require the following result. 

\noindent
{\bf Lemma 5.}
{\it If $J_0J_1J_2 \cdots J_{n-1}J_0$ is a cycle of length $n$, then there exists a periodic point $y$ of $f$ such that $f^i(y)$ belongs to $J_i$ for $i = 0, 1, \cdots, n-1$ and $f^n(y) = y$.}

\noindent
{\it Proof.}
Let $Q_n = J_0$.  Since $f(J_{n-1}) \supset J_0 = Q_n$, there is, by Lemma 4, a closed subinterval $Q_{n-1}$ of $J_{n-1}$ such that $f(Q_{n-1}) = Q_n = J_0$.  Continuing this process one by one, we obtain, for each $0 \le i \le n-1$, a closed subinterval $Q_i$ of $J_i$ such that $f(Q_i) = Q_{i+1}$.  Consequently, $f^i(Q_0) = Q_i$ for all $0 \le i \le n$.  In particular, $f^n(Q_0) = Q_n = J_0 \supset Q_0$.  By Lemma 2, there is a point $y$ in $Q_0 \subset J_0$ such that $f^n(y) = y$.  Since $y \in Q_0$, we also obtain that $f^i(y) \in f^i(Q_0) = Q_i \subset J_i$ for all $0 \le i \le  n-1$. 
\hfill\sq

We call $f$ turbulent if there exist two closed subintervals $I_0$ and $I_1$ of $I$ with at most one point in common such that $f(I_0) \cap f(I_1) \supp I_0 \cup I_1$ and call $f$ strictly turbulent if such $I_0$ and $I_1$ are disjoint.  

\noindent
{\bf Lemma 6.}
If there exist a point $c$, a fixed point $z$ of $f$, and an integer $k \ge 2$ such that $f(c) < c < z \le f^k(c)$ \, ($f(c) < c < z < f^k(c)$ respectively), then there exist a point $d$, a fixed point $\hat z$ of $f$ such that $f(d) < d < \hat z \le f^2(d)$ \, $(f(d) < d < \hat z < f^2(d)$ respectively) and $f$ has periodic points of all periods and is turbulent (strictly turbulent respectively).    

\noindent
{\it Proof.}
If $f(c) < c < z \le f^2(c)$, by taking $I_0 = [f(c), c]$ and $I_1 = [c, z]$, we obtain that $f$ is turbulent.  If $f(c) < c < z < f^2(c)$, then there is a point $r$ such that $f(c) < f(r) < c < r < z < f^2(r) \approx f^2(c)$.  By taking $I_0 = [f(r), c]$ and $I_1 = [r, z]$, we see that $f$ is strictly turbulent.  In either case, by applying Lemma 5 to the cycles $I_0(I_1)^{i-1}I_0$, $i \ge 1$, we obtain periodic points of all periods for $f$.  Now suppose the lemma is true respectively for $k = n \ge 2$ and let $k = n+1$.  If $c < f^2(c) < z$ and $f^2(c) < f^3(c)$, let $z^*$ be a fixed point of $f$ in $(c, f^2(c))$.  If $z < f^2(c)$, let $z^* = z$.  In either case, we have $f(c) < c < z^* < f^2(c)$ and we are done as before.  If $f^2(c) < f(c)$, let $u = f(c)$.  If $f(c) < f^2(c) < c < z$, let $u$ be a point in $(c, z)$ such that $f(u) = f^2(c)$.  If $f(c) < c < f^2(c) < z$ and $f^3(c) < f^2(c)$, let $u$ be a point in $(f^2(c), z)$ such that $f(u) = f^2(c)$.  In either case, we have $f(u) < u < z \le f^{n+1}(c) = f^n(u)$ (or $f(u) < u < z < f^{n+1}(c) = f^n(u)$ respectively) and the lemma follows by induction.
\hfill\sq

\section{A non-directed graph proof of (a), (b) and (c)}
The proof we present here is slightly different from the one in {\bf{\cite{du3}}}.  It is more direct and intuitive.  It is especially suitable for the calculus course.  When reading this proof, it is better to draw an accompanying graph.  It will help to understand why (a), (b) and (c) hold so naturally.  

Let $P$ be a period-$m$ orbit of $f$ with $m \ge 3$.  Let $v$ be a point such that $\min P \le f^2(v) < v < f(v) = b \in P$.  Let $z$ be a fixed point of $f$ in $[v, b]$.  Since $f^2(\min P) > \min P$ and $f^2(v) < v$, the point $y = \max \{ \min P \le x \le v : f^2(x) = x \}$ exists.  Furthermore, $f(x) > z$ on $[y, v]$ and $f^2(x) < x$ on $(y, v]$.  Therefore, $y$ is a period-2 point of $f$.  (a) is proved.

For the proofs of (b) and (c), assume that $m \ge 3$ is odd and note that $f(x) > z > x > f^2(x)$ on $(y, v]$.  Since $f^{m+2}(y) = f(y) > y$ and $f^{m+2}(v) = f^2(v) = f(b) < v$, the point $p_{m+2} = \min \{ y \le x \le v : f^{m+2}(x) = x \}$ exists.  Let $k$ denote the least period of $p_{m+2}$ with respect to $f$.  Then $k$ divides $m+2$ and so $k$ is odd.  Furthermore, $k > 1$ because $f$ has no fixed points in $(y, v)$.  If $k < m+2$, let $x_k = p_{m+2}$, then $x_k$ is a solution of the equation $f^k(x) = x$ in $(y, v)$.  Since $f^{k+2}(y) = f(y) > y$ and $f^{k+2}(x_k) = f^2(f^k(x_k)) = f^2(x_k) < x_k$, the equation $f^{k+2}(x) = x$ has a solution $x_{k+2}$ in $(y, x_k)$.  Inductively, for each $n \ge 1$, the equation $f^{k+2n}(x) = x$ has a solution $x_{k+2n}$ such that $y < \cdots < x_{k+4} < x_{k+2} < x_k < v$.  Consequently, the equation $f^{m+2}(x) = x$ has a solution $x_{m+2}$ such that $y < x_{m+2} < x_k = p_{m+2}$.  This contradicts the minimality of $p_{m+2}$.  So, $p_{m+2}$ is a period-$(m+2)$ point of $f$.  This proves (b).  

We now prove (c).  Let $z_0 = \min \{ x : v \le x \le z$, $f^2(x) = x \}$.  Then $f^2(x) < x$ and $f(x) > z$ when $y < x < z_0$.  If $f^2(x) < z_0$ whenever $\min I \le x < z_0$, then $f^2([\min I, z_0]) \subset [\min I, z_0]$, which contradicts the fact that $(f^2)^{(m+1)/2}(v) = b > z_0$.  So, the point $d = \max \{ x : \min I \le x \le y, f^2(x) = z_0 \}$ exists and $f(x) > z \ge z_0 > f^2(x)$ for all $x$ in $(d, y)$.  Therefore $f(x) > z \ge z_0 > f^2(x)$ whenever $d < x < z_0$.  Let $s = \min \{ f^2(x) : d \le x \le z_0 \}$.  If $s \ge d$, then $f^2([d, z_0]) \subset [d, z_0]$, which again contradicts the fact that $(f^2)^{(m+1)/2}(v) = b > z_0$.  Thus $s < d$.  Let $u= \min \{ d \le x \le z_0 : f^2(x) = d \}$.  Since $f^2(d) = z_0 > d$ and $f^2(u) = d < u$, the point $c_2 = \min \{ d \le x \le u : f^2(x) = x \}$ \, $( \le y)$ exists and $d < f^2(x) < z_0$ on $(d, c_2]$.  Let $w$ be a point in $(d, c_2)$ such that $f^2(w) = u$.  Then since $f^4(d) = z_0 > d$ and $f^4(w) = d < w$, the point $c_4 = \min \{ d \le x \le w : f^4(x) = x \}$ \, $( < c_2)$ exists and $d < f^4(x)$ $< z_0$ on $(d, c_4]$.  Inductively, for each $n \ge 1$, let $c_{2n+2} = \min \{ d \le x \le c_{2n} : f^{2n+2}(x) = x \}$.  Then $d < \cdots < c_{2n+2} < c_{2n} < \cdots < c_4 < c_2 \le y$ and $d < (f^2)^k(x) < z_0$ on $(d, c_{2k}]$ for all $1 \le k \le n$.  Since $f(x) > z_0$ on $(d, z_0)$, we have $f^i(c_{2n}) < z_0 < f^j(c_{2n})$ for all even $2 \le i \le 2n$ and all odd $1 \le j \le 2n-1$.  So, each $c_{2n}$ is a period-$(2n)$ point of $f$.  Therefore, $f$ has periodic points of all even periods.  This establishes (c).  

\noindent
{\bf Remarks.}  (1) In the proof of (a), there are many ways to choose such point $v$.  One is: Let $a = \max \{ x \in P \, : \, f(x) > x \}$ and $b = \min \{ x \in P \, : \, a < x \}$.  Then $f(b) \le a < b \le f(a)$.  We can choose $v$ to be a point in $[a, b]$ such that $f(v) = b$.  The other is: Let $b$ be the point in $P$ such that $f(b) = \min P$ and let $v$ be a point in $[\min P, b]$ such that $f(v) = b$ (see section 10).

(2) The arguments in the proofs of (a) and (b) can be used to give a simpler proof of the main result of Block in {\bf{\cite{bl0}}} on the stability of periodic orbits in Sharkovsky's theorem.  Indeed, assume that $f$ has a period-$2^n$ point.  Let $F = f^{2^{n-2}}$.  Then $F$ has a period-4 orbit $P$.  Arguing as in the proof of (a), there exists a point $v$ such that $\min P \le F^2(v) < v < F(v)$ and $F^2(\min P) > \min P$.  So, there is an open neighborhood $N$ of $f$ in $C^0(I, I)$ such that, for each $g$ in $N$, the map $G = g^{2^{n-2}}$ satisfies $G^2(v) < v < G(v)$ and $G^2(\min P) > \min P$.  Thus, the point $y = \max \{ \min P \le x \le v : G^2(x) = x \}$ is a period-2 point of $G$.  By Lemma 3, $y$ is a period-$2^{n-1}$ point of $g$.  On the other hand, assume that $f$ has a period-$2^i \cdot m$ point with $m \ge 3$ and odd and $i \ge 0$.  Let $F = f^{2^i}$.  Then $F$ has period-$m$ points.  Arguing as in the proof of (b), there exist a period-2 point $y$ of $F$ and a point $v$ with $y < v$ such that $F^{m+2}(y) > y$ and $F^{m+2}(v) < v$ and $F$ has no fixed points in $[y, v]$.  So, there is an open neighborhood $N$ of $f$ in $C^0(I, I)$ such that, for each $g$ in $N$, the map $G = g^{2^i}$ satisfies that $G^{m+2}(y) > y$, $G^{m+2}(v) < v$, and $G$ has no fixed points in $[y, v]$.  By Sharkovsky's theorem, $g$ has period-$(2^i \cdot (m+2))$ points.

(3) In proving (c), we implicitly use the fact that $f^2(v) < v < z < f(v) = b = f^{m+1}(v) = (f^2)^{(m+1)/2}(v)$.  Surprisingly, these inequalities imply, by Lemma 6, the existence of periodic points of all periods for $f^2$.  However, the existence of periodic points of all periods for $f^2$ does not automatically guarantee the existence of periodic points of all {\it even} periods for $f$.  It only guarantees the existence of periodic points of $f$ with least period $2k$ for each even $k \ge 2$ and least period $\ell$ or $2\ell$ for each odd $\ell \ge 1$.  We need a little more work to ensure the existence of periodic points of all even periods for $f$ as we have done in the previous section and in {\bf{\cite{du3}}}.  In the next two sections, we present two different strategies by doing some suitable "surgery" to the map $f$ to achieve our goal.  These two proofs of (c) can be combined with the proofs of (a) and (b) in the previous section to give more non-directed graph proofs of (a), (b) and (c).  

In the following two sections, we shall let $m \ge 3$ be odd and let $a$ and $b$ be any two {\it ajacent} points in $P$ (so that $(a, b) \cap P = \emptyset$) such that $f(b) \le a < b \le f(a)$.  We also let $v$ be a point in $(a, b)$ such that $f(v) = b$.

\section{The first 'surgical' proof of (c)}
Let $z_1 \le z_2$ be the smallest and largest fixed points of $f$ in $[v, b]$ respectively.  Let $g$ be the continuous map on $I$ defined by $g(x) = \max \{ f(x), z_2 \}$ if $x \le z_1$; $g(x) = \min \{ f(x), z_1 \}$ if $x \ge z_2$; and $g(x) = -x + z_1 + z_2$ if $z_1 \le x \le z_2$.  Then $g([\min I, z_1]) \subset [z_2, \max I]$ and $g([z_2, \max I]) \subset [\min I, z_1]$ and $g^2(x) = x$ on $[z_1, z_2]$.  So, $g$ has no periodic points of any odd periods $\ge 3$.  On the other hand, it follows from the definition of $g$ that if the iterates $t, g(t), g^2(t), \cdots, g^k(t)$ or the iterates $t, f(t), f^2(t), \cdots, f^k(t)$ of some point $t$ are jumping alternatively between $[\min I, z_1)$ and $(z_2, \max I]$, then $g^i(t) = f^i(t)$ for all $1 \le i \le k$.  

Since $m \ge 3$ is odd, for some $1 \le i \le m-1$, both $f^i(b)$ and $f^{i+1}(b)$ lie on the same side of $[z_1, z_2]$.  Let $k$ be the {\it smallest} such $i$.  Then the iterates $b, f(b), f^2(b), \cdots, f^k(b)$ are jumping alternatively between $[\min I, z_1)$ and $(z_2, \max I]$ and since $f(v) = b$, so are the iterates $v, f(v), f^2(v), \cdots, f^{k+1}(v)$.  Consequently, $g^i(v) = f^i(v)$ for all $0 \le i \le k+1$.  If $k$ is odd, then $f^k(b) < z_1$ and $f^{k+1}(b) < z_1$ and so $g^{k+1}(b) = z_2$ and $g^2(v) = f^2(v) = f(b) < v < z_1 = g^{k+3}(v)$.  If $k$ is even, then $f^k(b) > z_2$ and $f^{k+1}(b) > z_2$ and so $g^{k+1}(b) = z_1$ and $g^2(v) = f^2(v) = f(b) < v < z_1 = g^{k+2}(v)$.  In either case, we have $g^2(v) < v < z_1 = (g^2)^n(v)$ for some $n \ge 2$.  By Lemma 6, $g^2$ has periodic points of all periods.  So, $g^2$ has period-$j$ points for all odd $j \ge 3$.  By Lemma 3, $g$ has either period-$j$ points or period-$(2j)$ points for any odd $j \ge 3$.  Since $g$ has no periodic points of any odd periods $\ge 3$, $g$ has period-$(2j)$ points which are also period-$(2j)$ points of $f$ for all odd $j \ge 3$.  This proves (c).

\section{The second 'surgical' proof of (c)}
We first consider the special case when $P$ is a period-$m$ orbit of $f$ such that $(\min P, \max P)$ contains no period-$m$ orbits of $f$.  By (b), $[\min P, \max P]$ contains a period-$(m+2)$ orbit $Q$ of $f$.  Let $g$ be the continuous map from $I$ into itself defined by $g(x) = \min Q$ if $f(x) \le \min Q$; $g(x) = \max Q$ if $f(x) \ge \max Q$; and $g(x) = f(x)$ elsewhere.  Then $g$ has no period-$m$ points and by (b) $g$ has no period-$j$ points for any odd $3 \le j \le m-2$.  Since $g$ has the period-$(m+2)$ orbit $Q$, there exist a fixed point $z$ of $g$ and a point $v$ such that $g(v) \in Q$ and  $g^2(v) < v < z < g(v) = (g^2)^{(m+3)/2}(v)$.  By Lemma 6, $g^2$ has periodic points of all periods.  Since $g$ has no period-$j$ points for any odd $3 \le j \le m-2$, $g$ has period-$(2j)$ points for all odd $3 \le j \le m$ which are also periodic points of $f$ with the same periods.  Thus, (c) is proved in this special case.  The rest of the proof is to reduce the general case to this special one.

For the general case, we may assume that $I = [\min P, \max P]$ and $v$ is the {\it largest} point in $[\min P, b]$ such that $f(v) = b$.  By Lemma 6, we may also assume that $f$ has no fixed points in $[\min P, v] \cup [b, \max P]$.  Let $z_1 \le z_2$ be the smallest and largest fixed points of $f$ in $[v, b]$ respectively.  Without loss of generality, we may assume that $f(x) = x$ for all $z_1 \le x \le z_2$.  If $f$ has a period-2 point $\hat y$ in $[v, z_1]$, then $z_2 < f(\hat y) < b$.  In this case, let $t = \hat y$ and $u = f(\hat y)$ and let $h$ be the continuous map from $I$ into itself defined by $h(x) = f(x)$ for $x$ not in $[t, u]$ and $h(x) = -x + t + u$ for $x$ in $[t, u]$.  Otherwise, let $t = v$ and $u = z_2$ and let $h(x) = f(x)$ for $x$ in $I$.  Then it is clear that if $S$ is a period-$k$ orbit of $h$ with $k \ge 3$, then it is also a period-$k$ orbit of $f$.  For each period-$k$ orbit $S$ of $h$ with $k \ge 3$, since $[\min S, \max S]$ must contain a fixed point and a period-2 point of $h$, we have $[\min S, \max S] \supset [t, u]$.  Let $w = \min \{ \max C : C$ is a period-$m$ orbit of $h \}$.  Then by Lemma 1 $w$ is a periodic point of $h$ whose least period divides $m$.  It is obvious that $w$ cannot be a fixed point of $h$.  If the least period of $w$ is $< m$, let $W$ denote the orbit of $w$ under $h$.  Then by (b), $(\min W, \max W)$ contains a period-$m$ orbit of $h$ which contradicts the minimality of $w$.  Therefore, $W$ is a period-$m$ orbit of $h$ such that $(\min W, \max W)$ contains no period-$m$ orbits of $h$.  It follows from what we have just proved above that $h$ has period-$(2j)$ points for all odd $3 \le j \le m$ which are also periodic points of $f$ with the same periods.  (c) is proved.

In the following, we present seven different (but some are related) directed graph proofs of (a), (b) and (c).  In some cases, we shall need the easy fact that if $f$ has period-3 points then $f$ has periodic points of all periods.

\section{The first directed graph proof of (a), (b) and (c)}
The proof we present here is different from the one in {\bf{\cite{du2}}}.  Let $P = \{x_i : 1 \le i \le m \}$, with $x_1 < x_2 < \cdots < x_m$, be a period-$m$ orbit of $f$ with $m \ge 3$.  Let $x_s = \max \{x \in P : f(x) > x \}$.  Then $f(x_s) \ge x_{s+1}$ and $f(x_{s+1}) \le x_s$.  So, $f$ has a fixed point $z$ in $[x_s, x_{s+1}]$.  If for all integers $i$ such that $1 \le i \le m-1$ the points $f^i(x_s)$ and $f^{i+1}(x_s)$ lie on opposite sides of $z$, then $f([x_1, x_s] \cap P) \subset [x_{s+1}, x_m] \cap P$ and $f([x_{s+1}, x_m] \cap P) \subset [x_1, x_s] \cap P$.  Since $f$ is one-to-one on $P$, this implies that $f([x_1, x_s] \cap P) = [x_{s+1}, x_m] \cap P$ and $f([x_{s+1}, x_m] \cap P) = [x_1, x_s] \cap P$.  In particular, $f([x_1, x_s]) \supset [x_{s+1}, x_m]$ and $f([x_{s+1}, x_m]) \supset [x_1, x_s]$.  By applying Lemma 5 to the cycle $[x_1, x_s][x_{s+1}, x_m][x_1, x_s]$, we obtain a period-2 point of $f$.  

In tne sequel, assume that for some $1 \le r \le m-1$ the points $f^r(x_s)$ and $f^{r+1}(x_s)$ lie on the same side of $z$ (this includes the case when $m$ is odd).  We may also assume that $r$ is the {\it smallest} such integer (and so the iterates $x_s, f(x_s), f^2(x_s), \cdots, f^r(x_s)$ are jumping around $z$ alternatively).  The following argument is negligible and is included for the interested readers only.  If $r$ is even, let $J_0$ be the smallest closed interval which contains $f^i(x_s)$ for all even $i$ in $[0, r]$ and let $J_1$ be the smallest closed interval which contains $x_s$ and $f^j(x_s)$ for all odd $j$ in $[1, r-1]$.  If $r$ is odd, let $J_0$ be the smallest closed interval which contains $f^i(x_s)$ for all even $i$ in $[0, r-1]$ and let $J_1$ be the smallest closed interval which contains $f^j(x_s)$ for all odd $j$ in $[1, r]$.  By applying Lemma 5 to the cycle $J_0J_1J_0$, we obtain a period-2 point of $f$.  This proves (a).

Now suppose both $f^r(x_s)$ and $f^{r+1}(x_s)$ lie on the right side of $z$ (if they both lie on the left side of $z$, the proof is similar).  Then $r \le m-2$ and since $f(x_{s+1}) \le x_s < x_{s+1}$, we have $f^r(x_s) > x_{s+1}$.  Let $k$ be the {\it smallest} positive integer such that $f^k(x_s) \ge f^r(x_s)$.  Then $1 \le k \le r$.  If $k = 1$, then for each $n \ge 2$, by appealing to Lemma 5 for the cycle $$[x_{s+1}, f^r(x_s)]([x_s, x_{s+1}])^{n-1}[x_{s+1}, f^r(x_s)]$$ of length $n$ (here $([x_s, x_{s+1}])^{n-1}$ represents $n-1$ copies of $[x_s, x_{s+1}]$), we see that $f$ has periodic points of all periods $\ge 2$.  So, suppose $k \ge 2$.  Since the iterates $x_s, f(x_s), f^2(x_s), \cdots$, $f^k(x_s)$ are jumping around $z$ alternatively, we actually have $k \ge 3$.  So, $f^{k-1}(x_s) < z < f^{k-2}(x_s) < f^r(x_s) \le f^k(x_s)$.  Let $W = [f^{k-1}(x_s), z]$, $V = [z, f^{k-2}(x_s)]$ and $U = [f^{k-2}(x_s), f^r(x_s)]$.  For each {\it{even}} integer $n \ge 2$, we apply Lemma 5 to the cycle $U(WV)^{(n-2)/2}WU$ of length $n$ to obtain periodic points of all even periods, including period-2, for $f$.  This proves (a) and (c).  

Finally, let $f^j(x_s)$ be the unique point in $\{ f^i(x_s) : 1 \le i \le k-2$ and $i$ is odd $\}$ which is closest to the point $f^r(x_s)$.  Then $1 \le j \le k-2$ and $x_{s+1} \le f^{k-2}(x_s) \le f^j(x_s) < f^r(x_s)$.  Let $L = [f^j(x_s), f^r(x_s)]$ and let $[a : b]$ denote the closed interval with $a$ and $b$ as endpoints.  For each integer $n \ge m-1$, we apply Lemma 5 to the cycle $$[x_s, z][f(x_s) : z] [f^2(x_s) : z] \cdots [f^{j-1}(x_s) : z][z, f^{k-2}(x_s)][f^{k-1}(x_s), z]
L([x_s, x_{s+1}])^{n-j-3}[x_s, z]$$ of length $n$ to obtain a period-$n$ point of $f$ (note that if $m \ge 4$ is even, this means that if $f$ doesn't swap $\{ x_1, x_2, \cdots, x_{m/2} \}$ and $\{ x_{(m/2)+1}, x_{(m/2)+2}, \cdots, x_m \}$ then $f$ has periodic points of odd period $m-1$).  This proves (b) and hence completes the proofs of (a), (b) and (c).

\section{The second directed graph proof of (a), (b) and (c)}
The proof we present here is a variant of that in {\bf{\cite{du2}}} except the proof of (a).  Let $P = \{x_i : 1 \le i \le m \}$, with $x_1 < x_2 < \cdots < x_m$, be a period-$m$ orbit of $f$ with $m \ge 3$.  Let $x_s = \max \{x \in P : f(x) > x \}$.  Then $f(x_s) \ge x_{s+1}$ and $f(x_{s+1}) \le x_s$.  So, $f$ has a fixed point $z$ in $[x_s, x_{s+1}]$.  If for all integers $i$ such that $1 \le i \le m-1$ and $i \ne s$ the points $f(x_i)$ and $f(x_{i+1})$ lie on the same side of $z$, then $f([x_1, x_s] \cap P) \subset [x_{s+1}, x_m] \cap P$ and $f([x_{s+1}, x_m] \cap P) \subset [x_1, x_s] \cap P$.  Since $f$ is one-to-one on $P$, we have $f([x_1, x_s] \cap P) = [x_{s+1}, x_m] \cap P$ and $f([x_{s+1}, x_m] \cap P) = [x_1, x_s] \cap P$.  In particular, we have $f([x_1, x_s]) \supset [x_{s+1}, x_m]$ and $f([x_{s+1}, x_m]) \supset [x_1, x_s]$.  Consequently, $f$ has a period-2 point.    

Now assume that for some integer $t$ such that $1 \le t \le m-1$ and $t \ne s$ the points $f(x_t)$ and $f(x_{t+1})$ lie on opposite sides of $z$ (this includes the case when $m$ is odd).  Then $f([x_t, x_{t+1}])$ $\supp [x_s, x_{s+1}]$.  For simplicity, suppose $x_t < x_s$.  If $x_{s+1} \le x_t$, the proof is similar.  Let $q$ be the {\it smallest} positive integer such that $f^q(x_s) \le x_t$.  Then $2 \le q \le m-1$.  For $i = 0, 1, \cdots, q-1$, let $J_i = [z : f^i(x_s)]$.  For each $n \ge m+1$, we appeal to Lemma 5 for the cycle $J_0J_1 \cdots J_{q-1}$ $[x_t, x_{t+1}]$ $([x_s, x_{s+1}])^{n-q-1}J_0$ of length $n$ to confirm the existence of a period-$n$ point of $f$.  This proves (b).

On the other hand, since $x_t < x_s$ and $f(x_s) \ge x_{s+1}$, we may also assume that $t$ is the largest integer in $[1, s-1]$ such that $f(x_t) \le x_s$.  So, $f(x_i) \ge x_{s+1}$ for all $t+1 \le i \le s$.  If $f(x_i) \ge x_{t+1}$ for all $s+1 \le i \le m$, then $f^n(x_s) \ge x_{t+1}$ for all $n \ge 1$, contradicting the fact that $f^j(x_s) = x_t$ for some $1 \le j \le m-1$.  So, there is a {\it smallest} integer $\ell$ with $s+1 \le \ell \le m$ such that $f(x_\ell) \le x_t$.  If $x_{t+1} \le f(x_i) \le x_{\ell-1}$ for all $t+1 \le i \le \ell -1$, then $f^n(x_s) \ge x_{t+1}$ for all $n \ge 1$, again contradicting the fact that $f^j(x_s) = x_t$ for some $1 \le j \le m-1$.  So, there is an integer $k$ with $t+1 \le k \le \ell -1$ such that $f(x_k) \ge x_\ell$.  If $t+1 \le k \le s$, let $U = [x_t, x_k]$, $V = [x_k, z]$ and $W = [z, x_\ell]$.  For each {\it{even}} integer $n \ge 2$, we apply Lemma 5 to the cycle $U(WV)^{(n-2)/2}WU$ of length $n$ to establish the existence of periodic points of all even periods, including period-2, for $f$.  If $s+1 < k \le \ell -1$, by applying Lemma 5 to the cycles $[x_{s+1}, x_k]([x_k, x_\ell])^n[x_{s+1}, x_k]$, $n \ge 1$, we obtain periodic points of all periods $\ge 2$.   This proves (a) and (c) and hence completes the proofs of (a), (b) and (c).

\section{The third directed graph proof of (a), (b) and (c)}
The proof we present here follows the line of the usual 'standard' proof {\bf{\cite{al, bl, bh, ho}}} which involves periodic orbits of some special types called \v Stefan cycles.  However, our argument is different and simpler.

Let $P$ be a period-$m$ orbit of $f$ with $m \ge 3$ and odd.  If there is a point $p \in P$ such that either $$f^{m-2}(p) < \cdots < f^3(p) < f(p) < p = f^m(p) < f^2(p) < f^4(p) < \cdots < f^{m-1}(p)$$ or $$f^{m-1}(p) < \cdots < f^4(p) < f^2(p) < p = f^m(p) < f(p) < f^3(p) < \cdots < f^{m-2}(p),$$ then we say that $P$ is a \v Stefan cycle of $f$ with least period $m$.  

Let $P = \{x_i : 1 \le i \le m \}$, with $x_1 < x_2 < \cdots < x_m$, be a period-$m$ orbit of $f$ with $m \ge 3$.  Let $x_s = \max \{x \in P : f(x) > x \}$.  Then $f(x_s) \ge x_{s+1}$ and $f(x_{s+1}) \le x_s$.  So, $f$ has a fixed point $z$ in $[x_s, x_{s+1}]$.  If for all integers $i$ such that $1 \le i \le m-1$ the points $f^i(x_s)$ and $f^{i+1}(x_s)$ lie on opposite sides of $z$, then $f([x_1, x_s] \cap P) \subset [x_{s+1}, x_m] \cap P$ and $f([x_{s+1}, x_m] \cap P) \subset [x_1, x_s] \cap P$.  Since $f$ is one-to-one on $P$, this implies that $f([x_1, x_s] \cap P) = [x_{s+1}, x_m] \cap P$ and $f([x_{s+1}, x_m] \cap P) = [x_1, x_s] \cap P$.  In particular, $f([x_1, x_s]) \supset [x_{s+1}, x_m]$ and $f([x_{s+1}, x_m]) \supset [x_1, x_s]$.  By applying Lemma 5 to the cycle $[x_1, x_s][x_{s+1}, x_m][x_1, x_s]$, we obtain a period-2 point of $f$.  

On the other hand, we assume that there is a {\it smallest} integer $1 \le r \le m-1$ such that the points $f^r(x_s)$ and $f^{r+1}(x_s)$ lie on the same side of $z$ (this includes the case when $m$ is odd).  We may also assume that they both lie on the left side of $z$.  Then $2 \le r \le m-1$.  Since $f$ maps each of $x_s$ and $x_{s+1}$ to the other side of $z$, $f^r(x_s) \notin \{ \, x_s, x_{s+1}\, \}$.  Let $J_i = [z : f^i(x_s)]$ for all $0 \le i \le r-1$ and $J_r = [f^r(x_s), x_s]$ if $f^r(x_s)< x_s$ and $J_r = [x_{s+1}, f^r(x_s)]$ if $x_{s+1} < f^r(x_s)$.  Let $k > m$ be an odd integer.  Then we apply Lemma 5 to the cycle $J_0J_1 \cdots J_r([x_s, x_{s+1}])^{k-r-1}J_0$ of odd length $k$ to obtain a periodic point of odd period $\ell$ with $\ell > 1$ and $\ell$ divides $k$.  Therefore, for the proof of the rest of (a), it suffices to prove the case when $m \ge 3$ is odd.  Consequently, for the proof of (a), (b) and (c), it suffices to prove them for odd $m \ge 3$.  

In the following, we show that if $m \ge 3$ is odd and $f$ has no period-$\ell$ points with $\ell$ odd and $1 < \ell < m$, then $P$ is a \v Stefan orbit and $f$ has periodic points of all {\it even} periods and all periods $> m$.  Indeed, if $m = 3$, the proof is easy.  So, suppose $m > 3$ and $f$ has no period-$\ell$ points with $\ell$ odd and $1 < \ell < m$.  Let $J_i = [z : f^i(x_s)]$ for all $0 \le i \le r-1$ and $J_r = [f^r(x_s), x_s]$ if $f^r(x_s)< x_s$ and $J_r = [x_{s+1}, f^r(x_s)]$ if $x_{s+1} < f^r(x_s)$.  If $1 \le r \le m-3$, then we apply Lemma 5 to the cycle $J_0J_1 \cdots J_r([x_s, x_{s+1}])^{m-3-r}J_0$ of odd length $m-2$ to obtain a period-$\ell$ point of $f$ with $\ell$ odd and $1 < \ell < m$.  This is a contradiction.  Therefore, for all $1 \le k \le m-3$, the points $f^k(x_s)$ and $f^{k+1}(x_s)$ lie on opposite sides of $z$.  This implies that $f^i(x_s) < z < f^j(x_s)$ for all even $0 \le i \le m-3$ and all odd $1 \le j \le m-2$.  In particular, $z < f^{m-2}(x_s)$.   We now have four cases to consider depending on the locations of $f^{m-1}(x_s)$:

Case 1.  $f^{m-1}(x_s) < z < f^{m-2}(x_s)$.  In this case, we actually have $f^i(x_s) < z < f^j(x_s)$ for all even $0 \le i \le m-1$ and all odd $1 \le j \le m-2$.  For $0 \le i \le m-1$, let $J_i = [z : f^i(x_s)]$.  If, for some $1 \le j < k \le m-1$ with $k-j$ even, we have $[z : f^k(x_s)] \subset [z : f^j(x_s)]$, then by applying Lemma 5 to the cycle $J_0J_1 \cdots J_{j-1}J_kJ_{k+1} \cdots$ $J_{m-2}$ $[f^{m-1}(x_s), x_s]J_0$ if $k < m-1$ or to the cycle $J_0J_1J_2 \cdots J_{j-1} [f^{m-1}(x_s), x_s]J_0$ if $k = m-1$ of odd length $m - k + j$, we obtain a periodic point of $f$ of odd period strictly between 1 and $m$.  This contradicts the assumption.  So, we must have     
$$
f^{m-1}(x_s) < \cdots < f^4(x_s) < f^2(x_s) < x_s < f(x_s) < f^3(x_s) < \cdots < f^{m-2}(x_s).
$$
\noindent
That is, the orbit of $x_s$ is a \v Stefan cycle.  Let $U = [f^{m-1}(x_s), f^{m-3}(x_s)]$, $V = [f^{m-3}(x_s), z]$ and $W = [z, f^{m-2}(x_s)]$.  By applying Lemma 5 to the cycles $U(WV)^iWU$, $i \ge 0$, we obtain periodic points of all even periods for $f$.  On the other hand, for each integer $n > m$, by appealing to Lemma 5 for the cycle $J_0J_1J_2 \cdots J_{m-2} [f^{m-1}(x_s), f^{m-3}(x_s)] ([x_s, x_{s+1}])^{n-m} J_0$ of length $n$, we obtain a period-$n$ point for $f$.     

Case 2.  $x_{s+1} < f^{m-2}(x_s) < f^{m-1}(x_s)$.  In this case, there is a fixed point $\hat z$ of $f$ in $[x_{s+1}, f^{m-2}(x_s)]$.  By applying Lemma 5 to the cycle $[\hat z, f^{m-2}(x_s)]$ $([f^{m-2}(x_s), f^{m-1}(x_s)])^2$ $[\hat z, f^{m-2}(x_s)]$, we obtain a period-3 point of $f$ which is a contradiction.

Case 3.  $x_{s+1} < f^{m-1}(x_s) < f^{m-2}(x_s)$.  In this case, $x_{s+1} = f^k(x_s)$ for some $1 \le k \le m-3$.  By appealing to Lemma 5 for the cycle $[z, f^k(x_s)] [z : f^{k+1}(x_s)] [z : f^{k+2}(x_s)] \cdots [z : f^{m-3}(x_s)] [f^{m-1}(x_s), f^{m-2}(x_s)] ([x_s, x_{s+1}])^{k-1} [z, f^k(x_s)]$ of length $m-2$, we obtain a period-$\ell$ point of $f$ with $\ell$ odd and $1 < \ell < m$.  This contradicts the assumption.

Case 4.  $x_{s+1} = f^{m-1}(x_s) < f^{m-2}(x_s)$.  In this case, we have $f^i(x_s) < z < x_{s+1} = f^{m-1}(x_s) < f^j(x_s)$ for all even $0 \le i \le m-3$ and all odd $1 \le j \le m-2$.  Since $x_{s+1} = f^{m-1}(x_s)$, we have $f(x_{s+1}) = x_s$.  Thus, $x_s = f(x_{s+1})$.  By pluging this in the above inequalities, we obtain that $f^j(x_{s+1}) < z < f^i(x_{s+1})$ for all odd $1 \le j \le m-2$ and all even $0 \le i \le m-1$.  This is a symmetric copy of Case 1.  Therefore, $P$ is a \v Stefan cycle and $f$ has periodic points of all even periods and all periods $> m$.    

\noindent
This shows that $P$ is a \v Stefan cycle and $f$ has periodic points of all even periods and all periods $> m$.  

Now suppose $f$ has a periodic point of odd period $n \ge 3$.  Then there is an odd integer $1 < m \le n$ such that $f$ has a period-$m$ point and no periodic points of odd periods $\ell$ with $1 < \ell < m$.  It follows from what we have just proved that $f$ has periodic points of all even periods and for each $j \ge m+1$, in particular, for each $j \ge n+1$, $f$ has a period-$j$ point. This completes the proofs of (a), (b) and (c) simultaneously.
 
\section{The fourth directed graph proof of (a), (b) and (c)}
The proof we present here follows the line of the usual 'standard' proof {\bf{\cite{al, bl, bh, ho}}}.  However, our argument is different and doesn't need to know the structure of \v Stefan cycles beforehand.  

Let $P = \{x_i : 1 \le i \le m \}$, with $x_1 < x_2 < \cdots < x_m$, be a period-$m$ orbit of $f$ with $m \ge 3$.  Let $x_s = \max \{x \in P : f(x) > x \}$.  Then $f(x_s) \ge x_{s+1}$ and $f(x_{s+1}) \le x_s$.  So, $f$ has a fixed point $z$ in $[x_s, x_{s+1}]$.  If for all integers $i$ such that $1 \le i \le m-1$ and $i \ne s$ the points $f(x_i)$ and $f(x_{i+1})$ lie on the same side of $z$, then $f([x_1, x_s] \cap P) \subset [x_{s+1}, x_m] \cap P$ and $f([x_{s+1}, x_m] \cap P) \subset [x_1, x_s] \cap P$.  Since $f$ is one-to-one on $P$, we actually have $f([x_1, x_s] \cap P) = [x_{s+1}, x_m] \cap P$ and $f([x_{s+1}, x_m] \cap P) = [x_1, x_s] \cap P$.  In particular, we have $f([x_1, x_s]) \supset [x_{s+1}, x_m]$ and $f([x_{s+1}, x_m]) \supset [x_1, x_s]$.  By Lemma 5, we obtain a period-2 point of $f$.  Now assume that there is an integer $1 \le t \le m-1$ such that $t \ne s$ and the points $f(x_t)$ and $f(x_{t+1})$ lie on opposite sides of $z$ (this includes the case when $m \ge 3$ is odd).  Let $q$ be the smallest positive integer such that $f^q(x_s) \le x_t$ if $x_t \le x_s$ or $f^q(x_s) \ge x_{t+1}$ if $x_t \ge x_{s+1}$.  Then $1 \le q \le m-1$.  By applying Lemma 5 to the cycle $[x_s, z] [z, f(x_s)] [z : f^2(x_s)] \cdots [z : f^{q-1}(x_s)] [x_t, x_{t+1}] ([x_s, x_{s+1}])^i [x_s, z]$, $i \ge 0$, we obtain that $f$ has periodic points of all periods $> m$.  In particular, when $m \ge 3$ is odd, this establishes (b) and when $m \ge 4$ is even, this shows that $f$ has periodic points of {\it odd} periods $> m$.  Therefore, for the proof of the rest of (a), it suffices to prove the case when $m \ge 3$ is odd which turns out to tb a very special case of (c).  Consequently, to complete the proof of (a) and (c), it suffices to prove just (c).  

In the following, we show that if $m \ge 3$ is odd and $f$ has no periodic points of odd period $\ell$ with $1 < \ell < m$, then $f$ has periodic points of all {\it even} periods and $P$ is a \v Stefan orbit (see also {\bf{\cite{du4}}}).  Indeed, if $m = 3$, the proof is easy.  So, suppose $m > 3$.  Since $m$ is odd, there is a {\it smallest} integer $1 \le r \le m-1$ such that the points $f^r(x_s)$ and $f^{r+1}(x_s)$ lie on the same side of $z$ and, for some integer $t$ such that $1 \le t \le m-1$ and $t \ne s$, the points $f(x_t)$ and $f(x_{t+1})$ lie on opposite sides of $z$ (and so $f([x_t, x_{t+1}])$ $\supp [x_s, x_{s+1}]$).  For simplicity, we assume that $x_t < x_s$.  (If $x_{s+1} \le x_t$, then we consider the iterates of the point $x_{s+1}$ instead of the point $x_s$ and the proof is similar).  Let $q$ be the {\it smallest} positive integer such that $f^q(x_s) \le x_t$.  Then $2 \le q \le m-1$.  Let $J_i = [z : f^i(x_s)]$ for all $0 \le i \le q-1$.  If $2 \le q \le m-3$, we apply Lemma 5 to the cycle $J_0J_1J_2 \cdots J_{q-1}[x_t, x_{t+1}]([x_s, x_{s+1}])^{m-2-q-1}J_0$ of length $m-2$ to obtain a period-$\ell$ point with $\ell$ odd and $1 < \ell < m$.  This is a contradiction.  So, either $q = m-1$ or $q = m-2$.      

If $1 \le r \le m-3$, then since $f$ maps each of $x_s$ and $x_{s+1}$ to the other side of $z$, $f^r(x_s) \notin \{ \, x_s, x_{s+1}\, \}$.  Let $J_i = [z : f^i(x_s)]$ for all $0 \le i \le r-1$ and $J_r = [f^r(x_s), x_s]$ if $f^r(x_s)< x_s$ and $J_r = [x_{s+1}, f^r(x_s)]$ if $x_{s+1} < f^r(x_s)$.  Then we apply Lemma 5 to the cycle $J_0J_1 \cdots J_r([x_s, x_{s+1}])^{m-2-r-1}J_0$ of odd length $m-2$ to obtain a periodic point of odd period $\ell$ with $1 < \ell < m$.  This is a contradiction.  (Note that this argument can also be used to show that if $m \ge 4$ is even and $f$ has no periodic points of odd period $\ell$ with $1 < \ell < m$, then the iterates $x_s, f(x_s), \cdots f^{m-1}(x_s)$ are jumping alternately around a fixed point of $f$).       Therefore, $z < f^{m-2}(x_s)$ and hence $q = m-1$ and $f^{m-1}(x_s) = x_t < f^i(x_s) < x_s < z < f^j(x_s)$ for all even $0 \le i \le m-3$ and all odd $1 \le j \le m-2$.  Let $U = [f^{m-1}(x_s), f^{m-3}(x_s)]$, $V = [f^{m-3}(x_s), z]$ and $W = [z, f^{m-2}(x_s)]$.  By applying Lemma 5 to the cycles $U(WV)^iWU$, $i \ge 0$, we obtain periodic points of all even periods for $f$.

For the interested readers, we now prove that such orbit $P$ is actually a \v Stefan orbit.  However, note that, as shown above, we don't need this fact to prove (b) and (c).  For $0 \le i \le m-1$, let $J_i = [z : f^i(x_s)]$.  If, for some $1 \le j < k \le m-1$ with $k-j$ even, we have $[z : f^k(x_s)] \sub [z : f^j(x_s)]$, then by applying Lemma 5 to the cycle $J_0J_1 \cdots J_{j-1}J_kJ_{k+1} \cdots$ $J_{m-2}$ $[x_t, x_{t+1}]J_0$ of odd length $m - k + j < m$, we obtain a periodic point of $f$ of odd period strictly between 1 and $m$.  This contradicts the assumption.  So, we must have     
$$
f^{m-1}(x_s) < \cdots < f^4(x_s) < f^2(x_s) < x_s < f(x_s) < f^3(x_s) < \cdots < f^{m-2}(x_s).
$$
\noindent
That is, the orbit of $x_s$ is a \v Stefan cycle.

Here is the second different but related proof of the fact that $f$ has periodic points of all even periods.  For this proof, we may assume that, since $f(x_s) > z$, $t$ is the {\it largest} integer in $[1, s-1]$ such that $f(x_t) \le x_s$.  Consequently, we obtain that $f(x_i) > z$ for all $t+1 \le i \le s$ and $f(x_t) \le x_s$ and keep in mind that either $q = m-1$ or $q = m-2$.  

If $x_t = f^{m-2}(x_s)$, then since $f(x_t) \le x_s$, we have either $f^{m-1}(x_s) < f^{m-2}(x_s)$ or $f^{m-2}(x_s) < f^{m-1}(x_s) < x_s$.  If $f^{m-2}(x_s) < x_t$, then again since $f(x_t) \le x_s$, we must have $f^{m-2}(x_s) < f^{m-1}(x_s) = x_t < x_s$.  Therefore, when $q = m-2$, we have the following three cases to consider: 

Case 1. $x_t = f^{m-2}(x_s) < f^{m-1}(x_s) = f(x_t) < x_s = f^m(x_s) < z$.  In this case, since $f(x_i) > z$ for all $t+1 \le i \le s$, we have a contradiction.

Case 2. $f^{m-1}(x_s) < f^{m-2}(x_s) = x_t$.  In this case, since we also have $f(x_s) > z > x_s$, $f$ has a fixed point $\hat z$ in $[x_t, x_s]$.  By applying Lemma 5 to the cycle $[f^{m-1}(x_s), x_t]$ $[x_t, \hat z]$ $[x_t, \hat z]$ $[f^{m-1}(x_s), x_t]$, we obtain a period-3 point of $f$ which contradicts the assumption.

Case 3. $f^{m-2}(x_s) < f^{m-1}(x_s) = x_t < x_s$.  In this case, since $x_{s+1} = f^k(x_s)$ for some $1 \le k \le m-3$, by applying Lemma 5 to the cycle $[z, f^k(x_s)] [z : f^{k+1}(x_s)] [z : f^{k+2}(x_s)] \cdots [z : f^{m-3}(x_s)] [x_t, x_{t+1}] ([x_s, x_{s+1}])^{k-1} [z, f^k(x_s)]$ of length $m-2$, we obtain a period-$\ell$ point of $f$ with $\ell$ odd and $1 < \ell < m$.  This is a contradiction.

\noindent
In either case, we obtain a contradiction.  Therefore, $q = m-1$ and $f^{m-1}(x_s) = x_t$.  If $x_t < f^{m-2}(x_s) < x_s$, then arguing as in Case 2 above, $f$ has a period-3 point which is a contradiction.  So, $x_{s+1} \le f^{m-2}(x_s)$. Now consider the locations of $f^{m-3}(x_s)$.  If $ x_{s+1} \le f^{m-2}(x_s) < f^{m-3}(x_s)$, then by applying Lemma 5 to the cycle $[x_s, z] [z, f(x_s)] [z : f^2(x_s)] \cdots [z : f^{m-4}(x_s)] [f^{m-2}(x_s), f^{m-3}(x_s)] [x_s, z]$ of length $m-2$, we obtain a period-$\ell$ point of $f$ with $\ell$ odd and $1 < \ell < m$.  This is a contradiction.  If $x_{s+1} < f^{m-3}(x_s) < f^{m-2}(x_s)$, then by appealing to Lemma 5 for the cycle $[z, f^{m-3}(x_s)]$ $[f^{m-3}(x_s), f^{m-2}(x_s)]$ $[f^{m-3}(x_s), f^{m-2}(x_s)]$ $[z, f^{m-3}(x_s)]$, we obtain a period-3 point of $f$ which again is a contradiction.  Therefore, we have $f^{m-1}(x_s) = x_t < x_{t+1} \le f^{m-3}(x_s) < z < f^{m-2}(x_s)$.  
Let $U = [f^{m-1}(x_s), f^{m-3}(x_s)]$, $V = [f^{m-3}(x_s), z]$ and $W = [z, f^{m-2}(x_s)]$.  By applying Lemma 5 to the cycles $U(WV)^iWU$, $i \ge 0$, we obtain periodic points of all even periods for $f$.  

Now suppose $f$ has a periodic point of odd period $n \ge 3$.  Then there is an odd integer $1 < m \le n$ such that $f$ has a period-$m$ point and no periodic points of odd periods $\ell$ with $1 < \ell < m$.  It follows from what we have just proved that $f$ has periodic points of all even periods. This completes the proofs of (a) and (c) simultaneously.  

\section{The fifth directed graph proof of (a), (b) and (c)}
The proof we present here is very much different from those in the previous sections and depends heavily on the choice of the particular point $\min P$ in the given period-$m$ orbit $P$.

Let $P = \{x_i : 1 \le i \le m \}$, with $x_1 < x_2 < \cdots < x_m$, be a period-$m$ orbit of $f$ with $m \ge 3$.  Let $x_s = \max \{x \in P : f(x) > x \}$.  Then $f(x_s) \ge x_{s+1}$ and $f(x_{s+1}) \le x_s$.  Let $z$ be a fixed point of $f$ in $[x_s, x_{s+1}]$.  Let $v$ be a point in $[x_s, z]$ such that $f(v) = x_{s+1}$.  If for some point $x$ in $[f^2(v), v]$, $f(x) \le z$, then we apply Lemma 5 to the cycle $[x, v][z, f(v)][x, v]$ to obtain a period-2 point of $f$.  If for some $x$ in $[f^2(v), v]$, $f(x) > z$ and $v < f^2(x)$, then we appeal to Lemma 5 for the cycle $[x, v][f(v):f(x)][x, v]$ to obtain a period-2 point of $f$.  If $f^2(v) \le f^2(x) \le v < z < f(x)$ for all $x$ in $[f^2(v), v] \cap P$, then since $f^2$ is ono-to-one on $P$, we have $f^2([f^2(v), x_s] \cap P) = [f^2(v), x_s] \cap P$.  So, $f^2([f^2(v), x_s]) \supset [f^2(v), x_s]$.  Therefore, $f$ has a period-2 point.  The remaining case is that for some $v_1$ in $[f^2(v), v] \cap P$, $f^2(v_1) < f^2(v) < v_1 \le x_s \le v < z < f(v_1)$.  Now we can repeat the above argument with $v$ replaced by $v_1$ to obtain either a period-2 point of $f$ or a point $v_2$ in $P$ such that $f^2(v_2) < f^2(v_1) < v_2 < v_1 < z < f(v_2)$.  Proceeding in this manner finitely many times, we either have a period-2 point of $f$ or arrive at a point $v_k$ in $P$ such that $\min P = f^2(v_k) < v_k < z < f(v_k)$.  However, if there is a point $v_k$ in $P$ such that $\min P = f^2(v_k) < v_k < z < f(v_k)$, then the above argument implies that $f$ has a period-2 point.  Therefore, in either case, $f$ has a period-2 point.  This establishes (a).  Note that in the above proof of (a), if we can find a point $v$ such that $\min P = f^2(v) < v < f(v)$ in the first place, then the above argument will immediately show that $f$ has a period-2 point.  This is what we do below.  

Let $P$ be a period-$m$ orbit of $f$ with $m \ge 3$ and let $b$ be a point in $P$ such that $f(b) = \min P$.  If $f(x) < b$ for all $\min P \le x < b$, then $(\min P \le) \, f(\min P) < b$.  So, $(\min P \le) \, f^2(\min P) = f(f(\min P)) < b$.  By induction, $(\min P \le) \, f^k(\min P) < b$ for each $k \ge 1$.  This contradicts the fact that $b = f^m(b) = f^{m-1}(f(b)) = f^{m-1}(\min P) < b$.  Therefore, there exists a point $v$ (not necessarily a point of $P$) in $(\min P, b)$ such that $f(v) = b$.  Let $z$ be a fixed point of $f$ in $(v, b)$.    

If for some point $x$ in $[\min P, v]$, $f(x) \le z$, then we consider the cycle $[x, v][z, b][x, v]$.  If for some $x$ in $[\min P, v]$, $f(x) > z$ and $v < f^2(x)$, then we consider the cycle $[x, v][b : f(x)][x, v]$.  If $f^2(x) \le v < z < f(x)$ for all $x$ in $[\min P, v]$, let $q = \max([\min P, v] \cap P) \,\, (\le v)$.  Then $f^2([\min P, q] \cap P) \subset [\min P, q] \cap P$.  Since $f^2$ is one-to-one on $P$, this implies that $f^2([\min P, q] \cap P) = [\min P, q] \cap P$ and so $f^2([\min P, q]) \supset [\min P, q]$.  In this case, we consider the cycle $[\min P, q]f([\min P, q])[\min P, q]$.  In either of the above three cases, by Lemma 5, we obtain a period-2 point $y$ of $f$ such that $y < v < f(y)$.  This proves (a).   

Here is another proof of (a).  If there is a fixed point $\hat z$ of $f$ in $[\min P, v]$, then by applying Lemma 5 to the cycles $[v, b]([\hat z, v])^{n-1}[v, b]$, $n \ge 1$, we obtain periodic points of all periods for $f$ and we are done.  So, assume that $f$ has no fixed points in $[\min P, v]$.  Since $f^2(\min P) \ge \min P$ and $f^2(v) = \min P < v$, there is a point $y$ in $[\min P, v]$ such that $f^2(y) = y$.  Since $f$ has no fixed points in $[\min P, v]$ and $f(v) = b > v$, $y$ is a period-2 point of $f$ such that $y < v < f(y)$.  (a) is proved.  

In the following, let $m \ge 3$ be odd.  The first proof of (c) can be seen as a directed graph version of the proof of (c) in section 3 with the choice of the particular point $\min P$.  Let $z_0 = \min \{ v \le x \le b : f^2(x) = x \}$.  Then $f$ has neither fixed points nor period-2 points in $[v, z_0)$.  Since $f(v) > v$ and $f^2(v) < v$, we have $f(x) > x$ and $f^2(x) < x < z_0$ on $[v, z_0)$.  If $f^2(x) < z_0$ for all $\min P \le x \le v$, then $f^2(x) < z_0$ for all $\min P \le x < z_0$.  In particular, $(\min P \le) \, f^2(\min P) < z_0$ and so $(\min P \le) \, f^2(f^2(\min P)) < z_0$.  By induction, $(\min P \le) \, (f^2)^k(\min P) < z_0$ for each $k \ge 1$.  Since $m \ge 3$ is odd, this contradicts the fact that $(f^2)^{(m-1)/2}(\min P) = b > z_0$.  Hence $\max \{ f^2(x) : \min P \le x \le v \} \ge z_0$.  Let $I_0 = [\min P, v]$ and $I_1 = [v, z_0]$.  Then $f^2(I_0) \cap f^2(I_1) \supset I_0 \cup I_1$.  For each $n \ge 1$, we apply Lemma 5 to the cycle $I_1(I_0)^nI_1$ (with respect to $f^2$) of length $n+1$ to obtain a point $w$ in $[v, z_0]$ such that $(f^2)^{n+1}(w) = w$ and  $(f^2)^i(w) \in I_0$ for all $1 \le i \le n$.  Thus $f^{2i}(w) < v < w < f(w)$ for all $1 \le i \le n$.  Consequently, $w$ is a period-$(2n+2)$ point of $f$.  Therefore, $f$ has periodic points of all even periods $\ge 4$.  This proves (c).

We now give the second proof of (c).  For this proof, we choose $v$, in the beginning of this section, to be the {\it largest} point in $[\min P, b]$ such that $f(v) = b$ and then let $z_0 = \min \{ v \le x \le b : f^2(x) = x \}$.  Then $v < z_0 \le f(z_0) < b$ and $f^2(x) < x < z_0$ on $[v, z_0)$.  We choose a sequence $x_{-1}, x_{-2}, x_{-3}, \cdots$ of points so that $v < x_{-2} < x_{-4} < \cdots < z_0 \le f(z_0) < \cdots < x_{-3} < x_{-1} < b$ and $f(x_{-k}) = x_{-k+1}$ for $k = 1, 2, \cdots$ with $x_0 = v$.  If $f^2(x) < z_0$ whenever $\min P \le x \le v$, then $f^2(x) < z_0$ for all $\min P \le x < z_0$.  In particular, $(\min P) \le f^{2k}(\min P) < z_0$ for all $k \ge 1$, contradicting the fact that $z_0 < b = (f^2)^{(m-1)/2}(\min P)$.  So, $\max \{ f^2(x) : \min P \le x \le v \} \ge z_0$.  For each {\it even} integer $n \ge 2$, we appeal to Lemma 5 for the cycle $$[\min P, v]f([\min P, v])[x_{-n+2}, x_{-n}][x_{-n+3} : x_{-n+1}] \cdots [x_0 : x_{-2}][x_{-1}, b][\min P, v]$$ of length $n+2$ to obtain periodic points of all even periods $\ge 4$ for $f$.  This also proves (c).

Finally we prove (b).  Let $m \ge 3$ be odd.  Let $j$ be the largest {\it even} integer in $[0, m-1]$ such that $f^j(\min P) < y$.  Then $0 \le j < m-1$ since $f^{m-1}(\min P) = b > y$.  For $0 \le i \le m-1$, let $J_i = [f^i(\min P) : f^i(y)]$.  If there is an odd integer $\ell$ with $j < \ell < m-1$ such that $f^\ell(\min P) < y$, let $k$ denote the smallest such $\ell$.  Then $k-j \le m-2$ and $y < f^i(\min P)$ for all $j < i < k$.  We apply Lemma 5 to the cycle $$J_jJ_{j+1}J_{j+2} \cdots J_{k-1} [v, f(y)]([v, b])^{n+j-k-1}J_j$$ of length $n$ for each $n \ge m+1$.  If $y < f^\ell(\min P)$ for all odd integers $j < \ell < m-1$, we apply Lemma 5 to the cycle $$J_jJ_{j+1}J_{j+2} \cdots J_{m-2}([v, b])^{n-m+j+1}J_j$$ of length $n$ for each $n \ge m+1$.  In either case, for each $n \ge m+1$, we obtain a point $p_n$ in $[\min P, y]$ such that $f^n(p_n) = p_n$ and $p_n < y < f^i(p_n)$ for all $1 \le i \le n-1$.  Therefore, in either case, $p_n$ is a period-$n$ point of $f$ for each $n \ge m+1$.  This establishes (b).

\section{The sixth directed graph proof of (a), (b) and (c)}
The proof we present here is different from that in the previous section.  Let $P$ be a period-$m$ orbit of $f$ with $m \ge 3$ and let $a$ and $b$ be points in $P$ such that $f(a) = \max P$ and $f(b) = \min P$.  If $b < a$, let $z$ be the smallest fixed point of $f$ in $(b, a)$.  Then it is easy to see that $\max \{ f(x) : \min P \le x \le b \} \ge z$.   By applying Lemma 5 to the cycles $[b, z]([\min P, b])^i [b, z]$, $i \ge 1$, we obtain that $f$ has periodic points of all periods $\ge 2$ and we are done.  So, in the sequel, we suppose $a < b$.  

Let $v$ be the unique point in $P$ such that $f(v) = b$.  If $b < v$, then by applying Lemma 5 to the cycle $[b, v] ([a, b])^i [b, v]$, $i \ge 1$, we obtain that $f$ has periodic points of all periods $\ge 2$.  If $\min P < v < b$ and there is a fixed point $z$ of $f$ in $[\min P, v]$, then we apply Lemma 5 to the cycles $[z, v] ([v, b])^i [z, v]$, $i \ge 1$, to obtain that $f$ has periodic points of all periods $\ge 2$.  So, in the sequel, we assume that $\min P < v < b$ and $f$ has no fixed points in $[\min P, v]$.  

Since $f^2$ is one-to-one on $P$, we have $f^2([\min P, v] \cap P) \supset [\min P, v] \cap P$.  In particular, $f^2([\min P, v]) \supset [\min P, v]$.  By Lemma 2, there is a point $y$ in $[\min P, v]$ such that $f^2(y) = y$.  Since $f$ has no fixed points in $[\min P, v]$, $y$ is a period-2 point of $f$ such that $\min P < y < v < f(y)$.  This establishes (a).  

We now give a proof of (b) which is different from those given in the previous sections and can be seen as a directed graph version of the proof of (b) in section 3 with the chioce of the particular point $\min P$.  Let $m \ge 3$ be odd.  Since $f^{m+2}(y) = f(y) > v$ and $f^{m+2}(v) = \min P < y$, we have $[f^{m+2}(y) : f^{m+2}(v)] \supset [y, v]$.  Hence, by applying Lemma 5 to the cycle $$[y, v][f(y) : f(v)][f^2(y) : f^2(v)] [f^3(y) : f^3(v)] \cdots [f^{m+1}(y) : f^{m+1}(v)][y, v]$$ of length $m+2$, we obtain a point $r$ in $[y, v]$ such that $f^{m+2}(r) = r$, $f(r) \in [f(y) : f(v)]$ and $f^2(r) \in [f^2(y) : f^2(v)] = [\min P, y]$.  If $r$ has least period $m+2$, we are done.  Otherwise, let $n \ge 3$ be the least period of $r$ under $f$.  Then $n$ is odd and $m+2 - n > (m+2)/2$.  Since $[f^n(y) : f^n(r)] = [r, f(y)] \supset [v, f(y)]$ and $f([v, b]) \supset [v, b]$, we apply Lemma 5 to the cycle $$[y, r][f(y) : f(r)][f^2(y) : f^2(r)] \cdots [f^{n-1}(y) : f^{n-1}(r)][v, f(y)]([v, b])^{m+1-n}[y, r]$$ of length $m+2$ to obtain a point $s$ in $[y, r]$ such that $f^{m+2}(s) = s$, $f^2(s) \in [f^2(y): f^2(r)]$, $f^n(s) \in [v, f(y)]$ and $f^j(s) \in [v, b]$ for all $n+1 \le  j \le m+1$.  Since $f^2(r) < y$, $f^2(s) \in [f^2(y):f^2(r)] = [f^2(r), y]$.  Thus, $f^2(s) < y < s < v < f^\ell(s)$ for all $n \le \ell \le m+1$.  Since there are at least $(m+1)-n+1$ $(> (m+2)/2)$ {\it consecutive} iterates of $s$ under $f$ lying to the right of $y$ and $f^2(s) < y$, the least period of $s$ under $f$ is $m+2$.  This proves (b).

On the other hand, let $m \ge 3$ be odd and let $z_0 = \min \{ v \le x \le b : f^2(x) = x \}$.  If $f^2(x) < z_0$ for all  $\min P \le x \le v$, then $\min P \le f^{2i}(\min P) < z_0$ for all $i \ge 1$, contradicting the fact that $f^{2k}(\min P) = b > z_0$ for some $k \ge 1$.  So, $\max \{ f^2(x) : \min P \le x \le v \} \ge z_0$.  By applying Lemma 5 to the cycle $[v, z_0] ([\min P, v])^j [v, z_0]$, $j \ge 1$ (with respect to $f^2$), we obtain periodic points of $f$ of all even periods $\ge 4$.  So, (c) is proved.

In the following, we present two more different proofs of $(c)$ which are based on the existence of period-6 points for $f$.  Indeed, since $\min P = f(b) \le a < b \le \max P = f(a)$, there exist points $a < v < w < b$ such that $f(v) = b$ and $f(w) = a$.  Thus, if $m \ge 3$ is odd, then $P$ is also a period-$m$ orbit of $f^2$ and $\min P = f^2(v) < v < w < f^2(w) = \max P$.  It follows from the above that $f^2$ has period-3 points.  So, $f$ has either a period-6 point or a period-3 point.  This implies that $f$ has period-6 points if $m$ is odd.

We now prove that $f$ has period-$(2m)$ points.  If $m = 3$, then by applying Lemma 5 to the cycle $[y, v] [b : f(y)] [\min P, y] [f(\min P), f(y)] [v, b] [v, b] [y, v]$ of length 6, we obtain a period-6 point of $f$.  If $m \ge 5$ and odd, let $J_i = [f^i(\min P) : f^i(y)], i = 0, 1, 2, \cdots, m-1$.  Then since $[f^m(\min P):f^m(y)] \supset [y, f(y)]$, we apply Lemma 5 to the cycle $$[y, v][b : f(y)]J_0J_1J_2 \cdots J_{m-1}([y, f(y)])^{m-2}[y, v]$$ of length $2m$ to obtain a point $t$ such that $f^2(t) < y < t = f^{2m}(t) < v < f(t)$ and $f^i(t) \in [y, f(y)]$ for all $m+2 \le i \le 2m-1$.  Consequently, these $m+1$ points $t, f(t), f^2(t), f^i(t), m+2 \le i \le 2m-1$ are distinct points in the orbit of $t$.  So, $t$ is a period-$(2m)$ point of $f$.  This, together with the existence of period-6 points for $f$ proved above, establishes (c).

The next proof of (c) is based on the fact that if $f$ has period-6 points then $f$ has periodic points of all {\it even} periods.  Indeed, let $Q$ be a period-6 orbit of $f$.  Let $u$ and $w$ be any two {\it adjacent} points of $Q$ such that $f(w) \le u < w \le f(u)$.  Let $z$ be a fixed point of $f$ in $(u, w)$.  If there is no integer $i$ with $0 \le i \le 4$ such that the points $f^i(u)$ and $f^{i+1}(u)$ lie on the same side of $z$, then $f^i(u) < z < f^j(u)$ for $i = 0, 2, 4$ and $j = 1, 3, 5$.  We may assume that $f^2(u) < f^4(u) < u$ (if $f^4(u) < f^2(u) < u$, the proof is similar).  By applying Lemma 5 to the cycles $[f^4(u), u] [f(u), f^5(u)] [f^4(u), u]$ and $([f^4(u), u][f(u) : f^5(u)])^i [f^2(u), f^4(u)] [f^3(u) : f^5(u)] [f^4(u), u]$, $i \ge 1$, we obtain periodic points of all even periods for $f$.  So, suppose there is a {\it smallest} integer $r$ with $1 \le r \le 4$ such that the points $f^r(u)$ and $f^{r+1}(u)$ lie on the same side of $z$.  If $r = 1$ or $3$, choose points $u < u_3 < u_1 < z < u_2 < u_4 < w$ such that $f^i(u_1) = u_{i+1}$, $i = 1,2,3$ and $f(u_4) = u$.  By applying Lemma 5 to the cycle $[u_1, z] [z, u_2] [u_3, z] [z, u_4] [u, z] [w, f(u)] [u_1, z]$ if $r = 1$ or to the cycle $[u_3, z] [z, u_4] [u, z] [z, f(u)] [f^2(u), z] [w, f^3(u)] [u_3, z]$ if $r = 3$, we obtain a period-6 point $c$ of $f$ such that $f^i(c) < z < f^j(c)$ for $i = 0,2,4$ and $j = 1,3,5$.  It follows from above that $f$ has periodic points of all even periods.  If $r = 2$ or $4$, the proof is similar.  (c) is proved.

\section{The seventh directed graph proof of (a), (b) and (c)}
Let $P = \{x_i : 1 \le i \le m \}$, with $x_1 < x_2 < \cdots < x_m$, be a period-$m$ orbit of $f$ with $m \ge 3$.  Let $x_s = \max \{x \in P : f(x) > x \}$.  Then $f(x_s) \ge x_{s+1}$ and $f(x_{s+1}) \le x_s$.  If $b$ is a point (not necessarily a point of P) in $[\min P, x_s]$ such that $f(b) = \min P$, let $z$ be the smallest fixed point of $f$ in $[b, x_s]$.  Then it is easy to see that $\max \{ f(x) : \min P \le x \le b \} \ge z$.  By applying Lemma 5 to the cycles $[b, z] ([\min P, b])^i [b, z]$, $i \ge 1$, we obtain that $f$ has periodic points of all periods $\ge 2$ and we are done.  Same result holds if there is a point $a$ in $[x_{s+1}, \max P]$ such that $f(a) = \max P$.  So, the remaining case is that if $a$ and $b$ are two points in $P$ such that $f(a) = \max P$ and $f(b) = \min P$ then $a < x_s < x_{s+1} < b$.  

It is clear that one side of $[x_s, x_{s+1}]$ contains at least as many points of $P$ as the other side.  {\it We may assume that it is the right side}, (if it is the left side, the proof is similar).  If the right side contains as many points of $P$ as the left side, then we have $f([x_{s+1}, x_m] \cap P) \supset [x_1, x_s] \cap P$ and $f([x_1, x_s] \cap P) \supset [x_{s+1}, x_m] \cap P$ and so, $f$ has period-2 points.  If the right side contains more points of $P$ than the left side (this includes the case when $m \ge 3$ is odd), then it is clear that $f([x_{s+1}, x_m]) \supset [x_1, x_{s+1}]$ and $f([x_1, x_{s+1}]) \supset f([a, x_{s+1}]) \supset [x_{s+1}, x_m]$.  So, $f$ has period-2 points.  This establishes (a).  

Now since $\min P = f(b) \le a < b \le f(a) = \max P$, there exist points $a < v < w < b$ such that $f(v) = b$ and $f(w) = a$.  Thus, if $m \ge 3$ is odd, then $P$ is also a period-$m$ orbit of $f^2$ and $\min P = f^2(v) < v < w < f^2(w) = \max P$.  It follows from the above that $f^2$ has period-3 points.  So, $f$ has either a period-6 point or a period-3 point.  This implies that $f$ has period-6 points if $m$ is odd.

On the other hand, if $m$ is odd, then since $f(x_{s+1}) \le x_s$ and the right side of $[x_s, x_{s+1}]$ contains more points of $P$ than the left side, there is a smallest integer $s+1 \le t < m$ such that $f(x_{t+1}) \ge x_{s+1}$ and hence $f([x_t, x_{t+1}]) \supset [x_s, x_{s+1}]$.  Let $1 \le q \le m-1$ be the smallest integer such that $f^q(x_s) \ge x_{t+1}$.  By applying Lemma 5 to the cycles $[x_s, z] [f(x_s) : z] [f^2(x_s) : z] \cdots [f^{q-1}(x_s) : z]  [x_t, x_{t+1}] ([x_s, x_{s+1}])^{i-1} [x_s, z]$, $i \ge 1$, we obtain that $f$ has period-$n$ points for all $n \ge m$.  In particular, $f$ has period-$(m+2)$ points and period-$(2m)$ points.  This proves (b) and, together with the above, proves (c) also.  

For the sake of completeness, we include a proof of Sharkovsky's theorem which is slightly different from those in {\bf{\cite{du2}}}, {\bf{\cite{du3}}}.  

\section{A proof of Sharkovsky's theorem}

If $f$ has period-$m$ points with $m \ge 3$ and odd, then by (b) $f$ has period-$(m+2)$ points and by (c) $f$ has period-$(2 \cdot 3)$ points.  If $f$ has period-$(2 \cdot m)$ points with $m \ge 3$ and odd, then by Lemma 3, $f^2$ has period-$m$ points.  By (b), $f^2$ has period-$(m+2)$ points, which implies by Lemma 3 that $f$ has either period-$(m+2)$ points or period-$(2 \cdot (m+2))$ points.  If $f$ has period-$(m+2)$ points, then according to (c) $f$ has period-$(2 \cdot (m+2))$ points.  In either case, $f$ has period-$(2 \cdot (m+2))$ points.  On the other hand, since $f^2$ has period-$m$ points, by (c) $f^2$ has period-$(2 \cdot 3)$ points, hence by Lemma 3, $f$ has period-$(2^2 \cdot 3)$ points.  Now if $f$ has period-$(2^k \cdot m)$ points with $m \ge 3$ and odd and if $k \ge 2$, then by Lemma 3, $f^{2^{k-1}}$ has period-$(2 \cdot m)$ points.  It follows from what we have just proved that $f^{2^{k-1}}$ has period-$(2 \cdot (m+2))$ points and period-$(2^2 \cdot 3)$ points.  So, by Lemma 3, $f$ has period-$(2^k \cdot (m+2))$ points and period-$(2^{k+1} \cdot 3)$ points.  Consequently, if $f$ has period-$(2^i \cdot m)$ points with $m \ge 3$ and odd and if $i \ge 0$, then by Lemma 3, $f^{2^i}$ has period-$m$ points.  For each $\ell \ge i$, by Lemma 3, $f^{2^\ell} = (f^{2^i})^{2^{\ell -i}}$ has period-$m$ points.  By (c), $f^{2^\ell}$ has period-6 points.  So, $f^{2^{\ell +1}}$ has period-3 points and hence has period-2 points. This implies that $f$ has period-$2^{\ell+2}$ points for all $\ell \ge i$.  Finally, if $f$ has period-$2^k$ points for some $k \ge 2$, then $f^{2^{k-2}}$ has period-4 points and so, by (a), has period-2 points.  Therefore, $f$ has period-$2^{k-1}$ points.  This proves (1).

For the proofs of (2) and (3), it suffices to assume that $I = [0, 1]$.  Let $T(x) = 1 - |2x - 1|$ be the tent map on $I$.  Then for each positive integer $n$ the equation $T^n(x) = x$ has exactly $2^n$ distinct solutions in $I$.  It follows that $T$ has finitely many period-$n$ orbits.  Among these period-$n$ orbits, let $P_n$ be one with the smallest diameter $\max P_n - \min P_n$ (see also {\bf[1}, pp.32-34{\bf ]}).  For any $x$ in $I$, let $T_n(x) = \min P_n$ if $T(x) \le \min P_n$, $T_n(x) = \max P_n$ if $T(x) \ge \max P_n$, and $T_n(x) = T(x)$ elsewhere.  It is then easy to see that $T_n$ has exactly one period-$n$ orbit (i.e., $P_n$) but has no period-$m$ orbit for any $m$ with $m \prec n$ in the Sharkovsky ordering.  Now let $Q_3$ be the unique period-3 orbit of $T$ of smallest diameter.  Then $[\min Q_3, \max Q_3]$ contains finitely many period-6 orbits of $T$.  If $Q_6$ is one of smallest diameter, then $[\min Q_6, \max Q_6]$ contains finitely many period-12 orbits of $T$.  We choose one, say $Q_{12}$, of smallest diameter and continue the process inductively.  Let $q_0 = \sup \{\min Q_{2^i \cdot 3} : i \ge 0 \}$ and $q_1 = \inf \{ \max Q_{2^i \cdot 3} : i \ge 0 \}$.  Let $T_{\infty}(x) = q_0$ if $T(x) \le q_0$, $T_{\infty}(x) = q_1$ if $T(x) \ge q_1$, and $T_{\infty}(x) = T(x)$ elsewhere.  Then it is easy to see that $T_{\infty}$ has a period-$2^i$ point for $i = 0, 1, 2, \ldots$ but has no periodic points of any other periods.  This completes the proof of Sharkovsky's theorem.

\end{document}